\documentclass{amsart}

\usepackage{amssymb,mathtools,amsthm,xcolor}

\usepackage[numbers,sort&compress]{natbib}
\usepackage[
    colorlinks=true,
    linkcolor=blue,
    citecolor=purple,
    urlcolor=cyan,
    hyperfootnotes=false
    ]{hyperref}
\usepackage{tikz-cd}
\usepackage{cleveref}
\usepackage[margin=1in]{geometry}

\numberwithin{equation}{section}

\theoremstyle{plain}
\newtheorem{theorem}{Theorem}[section]
\newtheorem{proposition}[theorem]{Proposition}
\newtheorem{lemma}[theorem]{Lemma}

\newtheorem{conjecture}[theorem]{Conjecture}

\theoremstyle{definition}
\newtheorem{definition}[theorem]{Definition}

\theoremstyle{remark}
\newtheorem{remark}[theorem]{Remark}

\crefname{theorem}{Theorem}{Theorems}
\crefname{proposition}{Proposition}{Propositions}
\crefname{lemma}{Lemma}{Lemmas}
\crefname{corollary}{Corollary}{Corollaries}
\crefname{conjecture}{Conjecture}{Conjectures}
\crefname{definition}{Definition}{Definitions}
\crefname{example}{Example}{Examples}
\crefname{remark}{Remark}{Remarks}

\DeclareMathOperator{\Hom}{Hom}
\DeclareMathOperator{\Pic}{Pic}
\DeclareMathOperator{\Sym}{Sym}

\DeclareMathOperator{\rank}{rank}

\newcommand{\CC}{\mathbb{C}}
\newcommand{\PP}{\mathbb{P}}
\newcommand{\ZZ}{\mathbb{Z}}

\newcommand{\calA}{\mathcal{A}}
\newcommand{\calB}{\mathcal{B}}
\newcommand{\calE}{\mathcal{E}}
\newcommand{\calF}{\mathcal{F}}
\newcommand{\calK}{\mathcal{K}}
\newcommand{\calL}{\mathcal{L}}
\newcommand{\calM}{\mathcal{M}}
\newcommand{\calO}{\mathcal{O}}
\newcommand{\calQ}{\mathcal{Q}}
\newcommand{\calU}{\mathcal{U}}
\newcommand{\calV}{\mathcal{V}}
\newcommand{\contract}{\phi}

\newcommand{\printbibliography}{%
  \bibliographystyle{plainnat}%
  \bibliography{references}%
}

\title{Fujita's Very Ampleness Conjecture for $\PP^{n-1}$-Bundles over Curves}
\author{Kevin Dao}
\begin{document}

\maketitle

\begin{abstract}
    Let $\calE$ be a vector bundle of rank $n\geq 2$ over a smooth projective curve $C$ over $\CC$. We prove that $\omega_{\PP_C(\calE)}\otimes\calL_1\otimes\dots\otimes\calL_k$ is very ample for arbitrary ample line bundles $\calL_i$ on $\PP_C(\calE)$ when $k\geq n+2$. 
\end{abstract}

\section{Introduction}

In this short note, we prove a mixed form of Fujita's very ampleness conjecture for projective bundles over smooth projective curves in characteristic zero. Taking $\calL_1=\dots=\calL_k$ in \Cref{theorem:main} establishes Fujita’s very ampleness conjecture for these projective bundles. Recall that Fujita's conjectures comes in two parts.

\begin{conjecture}\label{conjecture:Fujita}
For a smooth projective variety $X$ of dimension $n$ over $\CC$ and an ample line bundle $\calL$, 
\begin{enumerate}
    \item $\omega_X\otimes\calL^{\otimes k}$ is globally generated if $k\geq n+1$,
    \item $\omega_X\otimes \calL^{\otimes k}$ is very ample if $k\geq n+2$.
\end{enumerate}
\end{conjecture}

Part (i) of \Cref{conjecture:Fujita} is known as Fujita's freeness conjecture. A proof of the freeness conjecture was announced in October, 2026  in a preprint by OpenAI \cite{OAI:Fujitas-freeness-conjecture-September-23-2026}. Part (ii) of \Cref{conjecture:Fujita} is known as Fujita's very ampleness conjecture and there has been limited progress on it since it was first conjectured.

Known cases of Fujita's very ampleness conjecture \cite{Fuj87} include curves, surfaces by Reider \cite{Rei88}, abelian varieties \cite[Example~2.1.9]{Mur19}, and smooth projective toric varieties. Musta\c{t}\u{a} proves a strengthened numerical statement in the smooth toric setting \cite{Mus02}, and Payne establishes the conjecture for the case of singular toric varieties under some added hypotheses \cite{Pay06}. If one assumes the line bundle $\calL$ in the conjecture were also globally generated, then the conclusion holds true by an argument using Castelnuovo--Mumford regularity and Kodaira vanishing as discussed in \cite[Example~2.1.10]{Mur19}.

\medskip

For the purposes of our note, we point out that earlier related work on projective bundles over curves includes Butler's numerical criteria for normal generation \cite{Butler1994} and Park's results on higher syzygies and mixed adjoints \cite{Park06}. 
Our proof builds on the methods used by Chen, K\"uronya, Mustopa, and Stix in their study of convex Fujita numbers of projective bundles \cite{chen2024convexfujitanumbersprojective}. 
Their work establishes global generation of the same mixed adjoint bundles for $k\geq n+1$. 
For related background on convex Fujita numbers, see also \cite{chen2023convexfujitanumberskodairaenriques}. 

The principal additional arguments needed here to get the very ampleness conjecture in our setting are the two lemmas in \Cref{ss:twolemmas}. 
These lemmas relate failure to separate a length-two subscheme to a rank condition on a particular morphism, and then to a slope inequality. 
These arguments require crucially that the base is a smooth curve.

\subsection*{Acknowledgements} KD acknowledges that the material is based upon work supported by the National Science Foundation Graduate
Research Fellowship Program under Grant No. 2137424. Any opinions, findings, and conclusions or recommendations expressed in this material are those of the authors and do not necessarily reflect the views of the National Science Foundation. 

KD thanks Ivan Aidun for reading a draft version of this note and providing comments which improved the exposition. KD also thanks Jack J. Garzella for their comments and suggestions.

\subsection*{AI Disclosure} AI assistants were used to search for, explain, and clarify various statements from the literature. AI assistants were also used to debug \LaTeX~code, generate figures, and proofread arguments. AI assistants were used also used to check the literature to determine if the statement was new. The writing and ideas are due to the author.

\subsection{Main Theorem}

We work over $\CC$, the field of complex numbers.

\begin{theorem}\label{theorem:main}
    Let $C$ be a smooth projective connected curve over $\CC$, $\calE$ be a vector bundle of rank $n\geq 2$ on $C$, and $\pi:\PP_C(\calE)\to C$ be the associated projective bundle over $C$. Let $\calL_1,\dots,\calL_k$ be ample line bundles on $\PP_C(\calE)$. If $k\geq n+2$, then
    \[
    \omega_{\PP_C(\calE)}\otimes \calL_1\otimes\dots\otimes\calL_k
    \]
    is a very ample line bundle.
\end{theorem}

\begin{remark}
    Establishing \Cref{theorem:main} also establishes the statement with $\CC$ replaced by an algebraically closed characteristic zero field by a standard Lefschetz principle argument.
\end{remark}

\subsection{Notation and preliminaries} 

This section summarizes our notations, terminologies, and recalls some preliminaries from the literature. First off, we fix $C$ a smooth projective connected curve. We let $\calE$ be a rank-$n$ vector bundle on $C$, $X:=\PP_C(\calE)$ the associated projective bundle, and $\pi:X\to C$ the structure morphism.

First, we start by recalling some properties of projective bundles. For more, we refer the reader to \cite[\href{https://stacks.math.columbia.edu/tag/01OA}{Tag 01OA}]{stacks-project}.

\begin{proposition}
The projective bundle $\pi:X\to C$ is a smooth projective variety of dimension $n$.
The Picard group of $X$ is
\[
\Pic(X)\cong \ZZ[\calO_X(1)]\oplus \pi^*\Pic(C).
\]
Furthermore, for $m\geq 0$,
\begin{equation}\label{eq:bundle-formulas}
\begin{aligned}
\pi_*\calO_X(m)
&\cong \Sym^m\calE\\
\omega_X &\cong \calO_X(-n)\otimes
\pi^*(\omega_C\otimes\det\calE).
\end{aligned}
\end{equation}
\end{proposition}

\begin{proof}
See \cite[\S2.3 and Equation 3.1]{chen2024convexfujitanumbersprojective} for the statement regarding the Picard group and the canonical bundle formula. See also \cite[Chapter III Exercise 8.4]{Hartshorne1977} for the isomorphism $\pi_*\calO_X(m)\cong \Sym^m\calE$.
\end{proof}

\begin{proposition}\cite[\href{https://stacks.math.columbia.edu/tag/0E8T}{Tag 0E8T}]{stacks-project}\label{proposition:VeryAmplenessCriteria}
A line bundle $\calL$ on a smooth projective variety $Y$ is very ample if it separates length $\leq 2$ subschemes. The latter condition means that for every length $\leq 2$ subscheme $Z\subseteq Y$, the restriction morphism
\[
H^0(Y,\calL)\to H^0(Z,\calL\otimes \calO_Z)
\]
is surjective.
\end{proposition}

\begin{definition}
    A vector bundle $\calV$ on $C$ is \textbf{(semi)-stable} if for all nonzero sub-vector bundles $\calV'\subsetneq \calV$, one has that
    \[ 
    \frac{\deg\calV'}{\rank\calV'} (\le) < \frac{\deg\calV}{\rank\calV}. 
    \]
    The quantity $\frac{\deg\calV}{\rank\calV}$ is called the \textbf{slope} of $\calV$ and denoted $\mu(\calV)$.
\end{definition}

\begin{figure}
    \centering
\begin{tikzpicture}[scale=0.75,
    x=1cm,y=1cm,
    >=Stealth,
    every node/.style={font=\small},
    polygon/.style={very thick},
    guide/.style={gray!45,densely dashed}
]
    \coordinate (O) at (0,0);
    \coordinate (F1) at (1.8,2.6);
    \coordinate (A) at (3.5,3.6);
    \coordinate (B) at (5.1,4.05);
    \coordinate (Flast) at (6.7,4.15);
    \coordinate (V) at (9,3.3);

    \draw[->] (-0.3,0) -- (10,0) node[right] {$\rank$};
    \draw[->] (0,-0.3) -- (0,5.1) node[above] {$\deg$};

    \draw[guide] (F1) -- (1.8,0)
        node[below=4pt,text=black] {$\rank(\calF_1)$};
    \draw[guide] (Flast) -- (6.7,0)
        node[below=4pt,text=black] {$\rank(\calF_{\ell-1})$};
    \draw[guide] (V) -- (9,0)
        node[below=4pt,text=black] {$\rank(\calV)$};

    \draw[polygon,blue!65!black] (O) -- (F1)
        node[midway,above,sloped] {$\mu_{\max}(\calV)$};
    \draw[polygon] (F1) -- (A);
    \draw[polygon,densely dotted] (A) -- (B);
    \draw[polygon] (B) -- (Flast);
    \draw[polygon,red!65!black] (Flast) -- (V)
        node[midway,above,sloped] {$\mu_{\min}(\calV)$};

    \fill (O) circle (1.8pt) node[below left] {$0$};
    \fill (F1) circle (1.8pt) node[below right=3pt] {$\calF_1$};
    \fill (Flast) circle (1.8pt) node[above=5pt] {$\calF_{\ell-1}$};
    \fill (V) circle (1.8pt) node[right=5pt] {$\calF_\ell=\calV$};
\end{tikzpicture}
    \caption{The Newton polygon of $\calV$.}
    \label{fig:HN-filtration}
\end{figure}
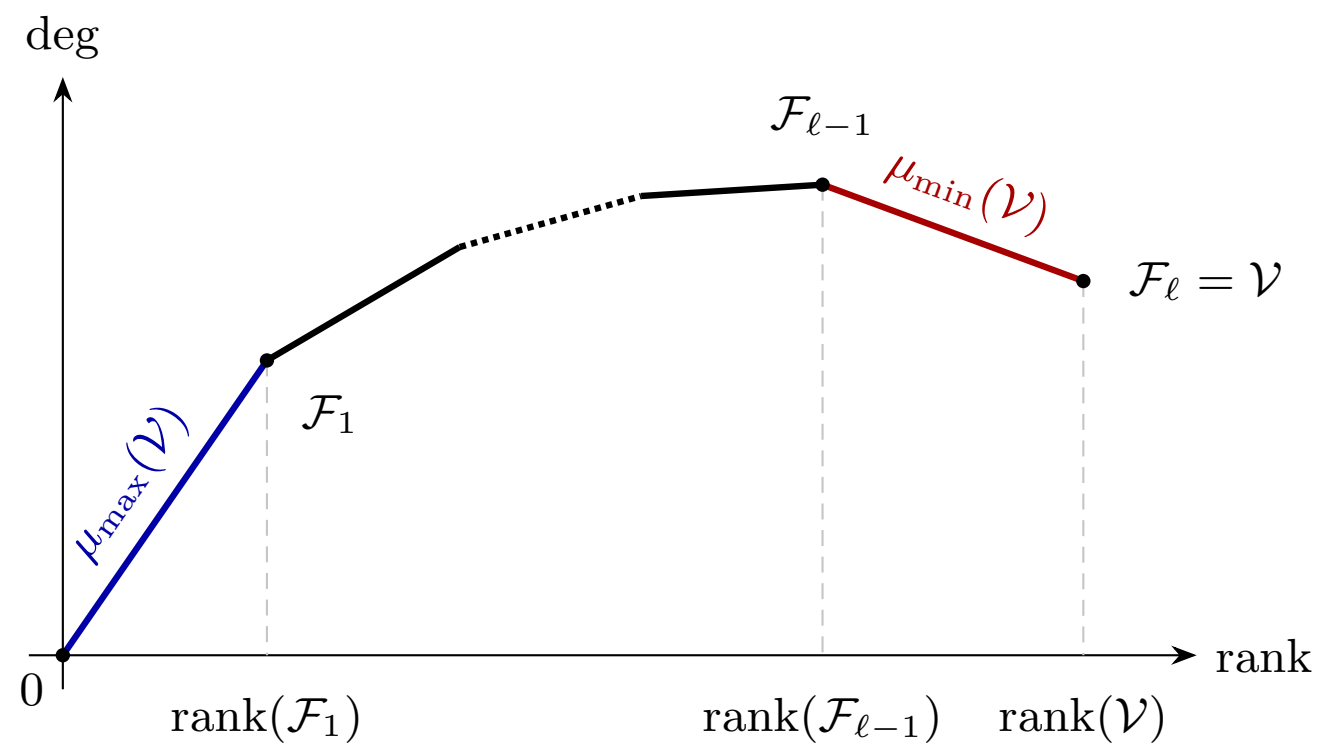

\begin{definition}
    For a vector bundle $\calV$ on $C$, its \textbf{Harder-Narasimhan filtration} is a filtration by subbundles
    \[
    0=:\calF_0\subsetneq \calF_1\subsetneq\calF_2\subsetneq\dots\subsetneq \calF_{\ell}=:\calV
    \]
    such that
    \begin{itemize}
        \item $\calF_i/\calF_{i-1}$ is semistable for all $1\leq i\leq \ell$,
        \item and the slopes $\mu_i=\mu(\calF_i/\calF_{i-1})$ are ordered in strictly decreasing order $\mu_1>\mu_2>\dots>\mu_\ell$.
    \end{itemize}
    It is a theorem of Harder and Narasimhan \cite{harder1975cohomology}  that such filtration exists, is unique, and that the numbers $\mu_1>\dots>\mu_\ell$ are uniquely determined.

    For this paper, we will use the notation $\mu_{\min}(\calV)$ and $\mu_{\max}(\calV)$ to denote the smallest and largest slopes occurring in the Harder-Narasimhan filtration. 
    
    We will not use it in this paper, but if one plots $(\rank\calF_i,\deg\calF_i)$, then one gets what is called the \textbf{Newton polygon} of the Harder-Narasimhan filtration (see \Cref{fig:HN-filtration}).
\end{definition}

\begin{proposition}\label{proposition:PropertiesOfHN}
The following properties regarding vector bundles $\calU,\calV$ on a smooth projective curve $C$ hold.
\begin{enumerate}
    \item If $\calV\twoheadrightarrow\calU$ is a nonzero quotient bundle, then $\mu_{\min}(\calV)\leq \mu_{\min}(\calU)$ and $\deg\calV\geq \rank(\calV)\mu_{\min}(\calV)$.
    \item For any $m\geq 0$, one has
    \begin{itemize}
        \item $\mu_{\max}(\calV^\vee)=-\mu_{\min}(\calV)$,
        \item $\mu_{\min}(\Sym^m\calV)=m\mu_{\min}(\calV)$,
        \item $\mu_{\max}(\calV\otimes \calU)=\mu_{\max}(\calV)+\mu_{\max}(\calU)$,
        \item $\mu_{\max}(\wedge^m\calV)\leq m\mu_{\max}(\calV)$ for $0\leq m\leq \rank(\calV)$ and,
        \item if $\calB$ is a line bundle of degree $b$ on $C$, then $\delta(\calB\otimes \calV)=\delta(\calV)+b$ for $\delta\in\{\mu,\mu_{\min},\mu_{\max}\}$.
    \end{itemize}
    \item If $H^0(C,\calV)\neq 0$, then $\mu_{\max}(\calV)\geq 0$.
    \item (Butler's Criterion) Let $\pi:X\to C$ with $X=\PP_C(\calE)$ be the projective bundle considered above. Let $\calM$ be a line bundle on $C$. Then, $\calO_X(a)\otimes \pi^*\calM$ is ample if and only if $a>0$ and $a\mu_{\min}(\calE)+\deg(\calM)>0$.
\end{enumerate}
\end{proposition}

\begin{proof}
The statements in (i), (ii), and (iii) are standard and we omit their proof. For (iv), see \cite[Proposition 3.5]{chen2024convexfujitanumbersprojective} which they attribute to \cite[Lemma 5.4]{Butler1994}. It is worth pointing out that (ii) critically uses the characteristic zero assumption.
\end{proof}

\begin{lemma}\label{lemma:ggSD}
Let $\calF$ be a vector bundle on $C$ with $\mu_{\min}(\calF)>1$. Then $\omega_C\otimes \calF$ is globally generated.
\end{lemma}

\begin{proof}
See \cite[Lemma 3.7]{chen2024convexfujitanumbersprojective} for (i).
\end{proof}

\section{Proof of Main Theorem}

For the rest of this note, we use $X:=\PP_C(\calE)$.

\begin{definition}
Let $m\geq 1$. Suppose we are given a nonzero morphism 
\[
q:\Sym^m\calE\to \calA
\]
with $\calA$ a line bundle on $C$. The \textbf{contraction} associated to $q$ is the morphism
\[
\phi_q:\calE\to (\Sym^{m-1}\calE)^\vee\otimes \calA
\]
given locally by $\phi_q(v)(v_2v_3\dots v_{m})=q(vv_2v_3\dots v_{m})$.
\end{definition}

\subsection{Two Lemmas}\label{ss:twolemmas}

\begin{lemma}\label{lemma:Estimate}
Let $m\geq 1$ and $q:\Sym^m\calE\to\calA$ be nonzero. Let $\calU:=\calE/\ker\contract_q$. Then,
\begin{enumerate}
    \item $\calU$ is a vector bundle on $C$ and,
    \item if $\contract_q$ has rank $\leq 1$ at some point of $C$, then 
    \[
    m\mu_{\min}(\calU)\leq \deg(\calA)-1+\frac{1}{\rank(\calU)}.
    \]
\end{enumerate}
\end{lemma}

\begin{proof}
Set $s:=\rank(\calU)$. First, $\calU$ is a vector bundle since it is the image of a morphism between vector bundles on a smooth curve. 

One can also check that $q$ vanishes when any of its $m$ inputs is in the kernel of $\phi_q$. Therefore, $q$ descends to $q_{\calU}:\Sym^m(\calU)\to \calA$ with $\phi_{q_{\calU}}$ generically injective. Furthermore, since $\calE\to \calU$ is surjective and splits locally, $\Sym^{m-1}\calE\twoheadrightarrow\Sym^{m-1}\calU$ and $(\Sym^{m-1}\calU)^\vee\hookrightarrow(\Sym^{m-1}\calE)^\vee$. Consequently, $\phi_q$ and $\phi_{q_\calU}$ have the same rank along each fibre.

By definition, $\phi_{q_\calU}$ is an element of $\Hom(\calU,(\Sym^{m-1}\calU)^\vee\otimes \calA)\cong H^0(C,\calU^\vee\otimes (\Sym^{m-1}\calU)^\vee\otimes \calA)$.
Therefore, its $s$th exterior power $w$ is a section of a vector bundle $T$ which we now compute. Taking the $s$th exterior power
of the generically injective morphism
\[
\phi_{q_{\calU}}:
\calU\longrightarrow
(\Sym^{m-1}\calU)^\vee\otimes\calA
\]
gives a nonzero morphism
\[
\bigwedge^s\phi_{q_{\calU}}:
\det\calU=\bigwedge^s\calU
\longrightarrow
\bigwedge^s\bigl((\Sym^{m-1}\calU)^\vee\otimes\calA\bigr).
\]
Thus
$\bigwedge^s\phi_{q_{\calU}}$ determines a nonzero section $w\in H^0(C,T)$, where
\begin{align*}
T :=
(\det\calU)^\vee\otimes
\bigwedge^s\bigl((\Sym^{m-1}\calU)^\vee\otimes\calA\bigr) \cong
(\det\calU)^\vee\otimes
\bigwedge^s\bigl((\Sym^{m-1}\calU)^\vee\bigr)
\otimes\calA^{\otimes s}.
\end{align*}
The last isomorphism uses the identity $\bigwedge^s(\calV\otimes\calA)
\cong
(\bigwedge^s\calV)\otimes\calA^{\otimes s}$ and this uses that $\calA$ is a line bundle.

Choose a point $p\in C$ where $\phi_{q_{\calU}}$ has
rank at most one, and let $\mathfrak m_p$ be the maximal
ideal of $\calO_{C,p}$.
Choose a local basis of $\calU$ whose last $s-1$ basis
vectors reduce to vectors in the kernel of the fibre map
at $p$.
In a matrix representing $\phi_{q_{\calU}}$, every entry
in these $s-1$ columns belongs to $\mathfrak m_p$.
Consequently, every $s\times s$ minor belongs to
$\mathfrak m_p^{s-1}$.
Thus the nonzero section $w$ vanishes to order at least
$s-1$ at $p$, so
\[
0\ne w\in H^0(C,T(-(s-1)p)).
\]

By \Cref{proposition:PropertiesOfHN} and $\deg(\calU)\geq s\mu_{\min}(\calU)$, we get
\begin{align*}
   s-1 & \leq \mu_{\max}(T) \leq s\deg(\calA)-\deg(\calU)-s(m-1)\mu_{\min}(\calU)\leq s(\deg(\calA)-m\mu_{\min}(\calU)).
\end{align*}
Dividing by $s$ and rewriting gives the desired inequality.
\end{proof}

\begin{lemma}\label{lemma:Obstruction}
Let $\calB\in\Pic(C)$. Define a vector bundle $\calF:=\Sym^m\calE\otimes \calB$ and a line bundle $\calQ:=\calO_X(m)\otimes\pi^*(\omega_C\otimes\calB)$ with $m\geq 1$. If $\mu_{\min}(\calF)>1$ and $\calQ$ does not separate a length-two subscheme $Z\subseteq X$, then there exists a degree-two effective divisor $D\subseteq C$ and a nonzero map $q:\Sym^m\calE\to \calB^\vee(D)$ whose contraction has rank at most one at a point of $D$.
\end{lemma}

\begin{proof}
By \Cref{lemma:ggSD}(i), $\omega_C\otimes\calF$ is globally generated. For a point $p\in C$, we write $X_p:=\pi^{-1}(p)$ for the fibre of $\pi$ over it. We have a commutative diagram
\[
\begin{tikzcd}
{H^0(X,\calQ)=H^0(X,\calO_X(m)\otimes \pi^*(\omega_C\otimes \calB))} \arrow[r] \arrow[d, "\cong" description] & {H^0(X_p,\calQ|_{X_p})} \arrow[d, "\cong" description] \\
{H^0(C,\omega_C\otimes \calF)=H^0(C,\Sym^m(\calE)\otimes \omega_C\otimes \calB)} \arrow[r, two heads]         & {H^0(p,\omega_C\otimes \calF|_p)}             
\end{tikzcd}
\]
where \eqref{eq:bundle-formulas} and the projection formula give the left arrow and the equalities there and global generation of $\omega_C\otimes\calF$ give us surjectivity of the bottom map. Together, we get that the restriction map $H^0(X,\calQ)\to H^0(X_p,\calQ|_{X_p})$ on top is surjective for every $p\in C$. So, if $\calQ$ does not separate $Z$, then $Z\not\subseteq X_p$ for any $p\in C$.

Consequently, $\pi|_Z$ identifies $Z$ with an effective divisor $D\subseteq C$ of degree two. For two distinct points, their images must be distinct. For a nonreduced $Z$, a local parameter on $C$ maps to a nonzero multiple of $\epsilon$ in $\calO_Z\cong\CC[\epsilon]/(\epsilon^2)$, giving the required isomorphism with a doubled point of $C$.

Let $\calK$ be the kernel of the natural evaluation map $\Sym^m\calE\to\pi_*(\calO_X(m)|_Z)$. This map is surjective because $Z$ is a section over $D$.

Moreover, $\calK$ is locally free and contains $(\Sym^m\calE)(-D)$.
Tensoring the evaluation sequence
\begin{equation}\label{eq:definingsequenceK}
0\to\calK\to\Sym^m\calE
\to\pi_*(\calO_X(m)|_Z)\to0.
\end{equation}
by the line bundle $\omega_C\otimes \calB$ gives the exact sequence
\[
0\longrightarrow\omega_C\otimes\calB\otimes\calK
\longrightarrow\omega_C\otimes\calF
\longrightarrow\pi_*(\calQ|_Z)
\longrightarrow0.
\]
Here, the cokernel is computed using the projection formula
\[
\omega_C\otimes \calB\otimes\pi_*(\calO_X(m)|_Z)\cong \pi_*(\pi^*(\omega_C\otimes \calB)\otimes \calO_X(m)|_Z)\cong \pi_*(\calQ|_Z).
\]
Since $\mu_{\min}(\calF)>1$, \Cref{proposition:PropertiesOfHN} (ii) implies $\mu_{\max}(\calF^\vee)<-1$ and so we have $H^0(C,\calF^\vee)=0$. By Serre duality, $H^1(C,\omega_C\otimes\calF)=0$. Since $\calQ$ does not separate $Z$,  the associated long exact sequence in cohomology and Serre duality gives
\[
0\ne
H^1(C,\omega_C\otimes\calB\otimes\calK)^\vee
\cong
\operatorname{Hom}(\calK,\calB^\vee).
\]
Choose a nonzero homomorphism $\psi:\calK\to\calB^\vee$. Restricting $\psi$ to $(\Sym^m\calE)(-D)$ also produces a nonzero map. Otherwise, $\psi$ factors through a nonzero map 
\[
\calK/(\Sym^m\calE)(-D)\to \calB^\vee
\]
but the former is a torsion sheaf supported on $D$ while $\calB^\vee$ is torsionfree and any such map is necessarily zero.

Restrict $\psi$ to $(\Sym^m\calE)(-D)$ and then twist by $D$ to get a nonzero map $q:\Sym^m\calE\to\calB^\vee(D)$. The restriction of $q$ to $\calK$,
\[
q|_{\calK}:\calK\to \calB^\vee(D)
\]
agrees with the composition
\[
\calK\stackrel{\psi}{\to}\calB^\vee\hookrightarrow\calB^\vee(D)
\]
since these maps are equal away from $D$.

Fix $p\in\operatorname{Supp}D$, and let $x$ be the unique closed point of $Z$ lying over $p$. For a vector bundle $\calV$ on $C$, write $\calV|_p:=\calV_p\otimes_{\calO_{C,p}}\kappa(p)$ for its fibre at $p$, where $\calV_p$ denotes its stalk. Similarly, $\calO_X(1)|_x$ denotes the one-dimensional fibre of $\calO_X(1)$ at $x$. 

The universal quotient map $\pi^*\calE\twoheadrightarrow\calO_X(1)$ induces a quotient map of vector spaces $\ell_x:\calE|_p\twoheadrightarrow\calO_X(1)|_x$. Since $Z\cong D$, tensoring \eqref{eq:definingsequenceK} with $\kappa(p)$ gives the right exact sequence
\[
\calK|_p
\longrightarrow
\Sym^m(\calE|_p)
\xrightarrow{\Sym^m(\ell_x)}
\bigl(\calO_X(1)|_x\bigr)^{\otimes m}
\longrightarrow0.
\]
Thus the image of the first arrow is
$\ker(\Sym^m(\ell_x))$.

The restriction of $q$ to $\calK$ factors through $\calB^\vee$ i.e.
\[
q|_{\calK}:\calK\to \calB^\vee\hookrightarrow\calB^\vee(D).
\]
This inclusion induces the zero map on fibres at $p$. Consequently, the induced fibre map $q|_p$ vanishes on $\ker(\Sym^m(\ell_x))$ and factors as
\[
\Sym^m(\calE|_p)
\xrightarrow{\Sym^m(\ell_x)}
\bigl(\calO_X(1)|_x\bigr)^{\otimes m}
\longrightarrow
\bigl(\calB^\vee(D)\bigr)|_p.
\]
For $v\in\ker(\ell_x)$ and $w\in\Sym^{m-1}(\calE|_p)$, this factorization gives $(q|_p)(vw)=0$. Hence the contraction of $q|_p$ vanishes on $\ker(\ell_x)$, a subspace of codimension one. So, the rank of $q$ at $p$ is at most one.
\end{proof}

\subsection{A Numerical Estimate}

\begin{proposition}\label{proposition:num_estimates}
    Let $\calL_1,\dots,\calL_k$ be ample line bundles on $X$ with $k\geq n+2$. Write $\calL_i=\calO_X(a_i)\otimes \pi^*\calM_i$ for $\calM_i\in\Pic(C)$ and $b_i:=\deg(\calM_i)$. Define $a:=\sum_{i=1}^k a_i$, $m:=a-n$, and $\calB:=\det(\calE)\otimes\bigotimes_{i=1}^k \calM_i$.

    Then,
    \begin{enumerate}
        \item $m\geq 2$ and $a_i\mu_{\min}(\calE)+b_i>0$ for each $1\leq i\leq k$,
        \item $\mu_{\min}(\Sym^m\calE\otimes\calB)>1$ and,
        \item for every nonzero quotient bundle $\calE\twoheadrightarrow\calU$, one has
        \[
        m\mu_{\min}(\calU)+\deg(\calB)\geq 1+\frac{2}{h}
        \]
        where $h$ is the rank of the final quotient in the Harder-Narasimhan filtration of $\calU$. 
    \end{enumerate}
\end{proposition}

\begin{proof}
By \Cref{proposition:PropertiesOfHN}(iv), each $a_i$ is positive and $a_i\mu_{\min}(\calE)+b_i>0$.
Thus $m=a-n\geq k-n\geq2$, proving (i).

The reduced denominator of $\mu_{\min}(\calE)$ is at most $n$. As $a_i$ and $b_i$ are integers, $a_i\mu_{\min}(\calE)+b_i\geq\frac{1}{n}$.
Consequently,
\begin{equation}\label{eq:lowerboundks}
\sum_{i=1}^k
\bigl(a_i\mu_{\min}(\calE)+b_i\bigr)
\geq \frac{k}{n}.
\end{equation}

We next prove (iii). Fix a nonzero quotient bundle $\calE\twoheadrightarrow\calU$, and let $h$ be the rank of $\calU$'s final Harder--Narasimhan quotient.
By \Cref{proposition:PropertiesOfHN}(i), $\mu_{\min}(\calU)\geq\mu_{\min}(\calE)$. Thus $a_i\mu_{\min}(\calU)+b_i>0$ for every $i$. Since $h\mu_{\min}(\calU)$ is an integer, we get $h\bigl(a_i\mu_{\min}(\calU)+b_i\bigr)\in\ZZ_{>0}$ and so
\[
\sum_{i=1}^k\bigl(a_i\mu_{\min}(\calU)+b_i\bigr)
\geq
\frac{k}{h}.
\]
Using $\deg\calB=\deg\calE+\sum_i b_i$, $\deg\calE\geq n\mu_{\min}(\calE)$, and $m\geq k-n$, we obtain two lower bounds:
\begin{align*}
m\mu_{\min}(\calU)+\deg\calB
&=
\sum_{i=1}^k \bigl(a_i\mu_{\min}(\calE)+b_i\bigr)
+\deg\calE-n\mu_{\min}(\calE)+m\bigl(\mu_{\min}(\calU)
-\mu_{\min}(\calE)\bigr)\\
&\geq
\frac{k}{n}
+(k-n)\bigl(\mu_{\min}(\calU)
-\mu_{\min}(\calE)\bigr),
\end{align*}
and
\begin{align*}
m\mu_{\min}(\calU)+\deg\calB
&=
\sum_{i=1}^k
\bigl(a_i\mu_{\min}(\calU)+b_i\bigr)
+\deg\calE-n\mu_{\min}(\calU)\geq
\frac{k}{h}
-n\bigl(\mu_{\min}(\calU)-\mu_{\min}(\calE)\bigr).
\end{align*}
Multiply the first inequality by $\frac{n}{k}$ and the second by $\frac{k-n}{k}$, then add. This gives
\[
m\mu_{\min}(\calU)+\deg\calB
\geq
\frac{n}{k}\frac{k}{n}
+\frac{k-n}{k}\frac{k}{h}
=
1+\frac{k-n}{h}
\geq1+\frac{2}{h}.
\]
and therefore proves (iii).

Finally, apply (iii) to the identity quotient
$\calE\twoheadrightarrow\calE$.
The formulas in \Cref{proposition:PropertiesOfHN}(ii) give
\[
\mu_{\min}(\Sym^m\calE\otimes\calB)
=
m\mu_{\min}(\calE)+\deg\calB
>1,
\]
proving (ii).
\end{proof}

\subsection{Proof of the Theorem}

\begin{proof}
Use $\calB$, $a$, and $m=a-n$ from \Cref{proposition:num_estimates}. By \eqref{eq:bundle-formulas},
\[
\omega_X\otimes\calL_1\otimes\dots\otimes\calL_k\cong \calQ
\]
where $\calQ:=\calO_X(a-n)\otimes\pi^*(\omega_C\otimes \calB)$. By \Cref{proposition:VeryAmplenessCriteria}, it suffices to prove the evaluation map $H^0(X,\calQ)\to H^0(Z,\calQ|_Z)$ is surjective for any length-two $Z\subseteq X$ subscheme. We show the dual of the cokernel of these maps must be zero.

Suppose that $\calQ$ fails to separate a length-two subscheme $Z$. By \Cref{lemma:Obstruction}, there is an effective divisor $D$ of degree two and a nonzero map $q:\operatorname{Sym}^m \mathcal{E} \rightarrow \mathcal{B}^{\vee}(D)$ whose contraction has rank at most one at a point of $D$.

Let $\calU:=\calE/\ker\phi_q$. Then \Cref{lemma:Estimate} together with $\deg(\calB^\vee(D))=2-\deg\calB$ implies
\begin{equation}\label{eq:mainthmeq1}
m\mu_{\min}(\calU)+\deg(\calB)\leq 1+\frac{1}{\rank(\calU)}.
\end{equation}
On the other hand, \Cref{proposition:num_estimates} says that
\begin{equation}\label{eq:mainthmeq2}
m\mu_{\min}(\calU)+\deg(\calB)\geq 1+\frac{2}{h}\geq 1+\frac{2}{\rank(\calU)}>1+\frac{1}{\rank(\calU)}
\end{equation}
where $h$ is the rank of the final quotient in the Harder-Narasimhan filtration of $\calU$. \Cref{eq:mainthmeq1} and \Cref{eq:mainthmeq2} are contradictory and we are done.
\end{proof}

\printbibliography
\end{document}